\documentstyle[12pt]{article}

\font\twlgot =eufm10 scaled \magstep1
\font\egtgot =eufm8
\font\sevgot =eufm7
\font\twlmsb =msbm10 scaled \magstep1
\font\egtmsb =msbm8
\font\sevmsb =msbm7

\newfam\gotfam

\textfont\gotfam\twlgot
\scriptfont\gotfam\egtgot
\scriptscriptfont\gotfam\sevgot

\newfam\msbfam
\textfont\msbfam\twlmsb
\scriptfont\msbfam\egtmsb
\scriptscriptfont\msbfam\sevmsb
\def\Bbb{\protect\pBbb}
\def\pBbb{\relax\ifmmode\expandafter\Bb\else\typeout{You cann't use
Bbb in text mode}\fi}
\def\Bb #1{{\fam\msbfam\relax#1}}

\def\thebibliography#1{\section*{References}\list
   {[\arabic{enumi}]}{\settowidth\labelwidth{#1}\leftmargin\labelwidth
     \advance\leftmargin\labelsep
     \usecounter{enumi}}
     \def\newblock{\hskip .11em plus .33em minus .07em}
     \sloppy\clubpenalty4000\widowpenalty4000
     \sfcode`\.=1000\relax}

\def\op#1{\mathop{\fam0 #1}\limits}

\newcommand{\beq}{\begin{equation}}
\newcommand{\eeq}{\end{equation}}
\newcommand{\ben}{\begin{eqnarray}}
\newcommand{\een}{\end{eqnarray}}
\newcommand{\be}{\begin{eqnarray*}}
\newcommand{\ee}{\end{eqnarray*}}
\newcommand{\bea}{\begin{eqalph}}
\newcommand{\eea}{\end{eqalph}}

\newcommand{\cH}{{\cal H}}

\newcommand{\cO}{{\cal O}}

\newcommand{\f}{\phi}
\newcommand{\om}{\omega}
\newcommand{\Om}{\Omega}

\newcommand{\g}{\gamma}
\newcommand{\G}{\Gamma}

\newcommand{\vf}{\varphi}

\newcommand{\si}{\sigma}

\newcommand{\w}{\wedge}
\newcommand{\wt}{\widetilde}

\newcommand{\dr}{\partial}

\newcounter{eqalph}
\newcounter{equationa}
\newcounter{theorem}
\newcounter{remark}
\newcounter{proposition}
\newcounter{lemma}
\newcounter{corollary}
\newcounter{definition}

\newenvironment{eqalph}{\stepcounter{equation}
\setcounter{equationa}{\value{equation}}
\setcounter{equation}{0}

\begin{eqnarray}}{\end{eqnarray}\setcounter{equation}{\value{equationa}}}

\def\theremark{\arabic{remark}}

\def\thedefinition{\arabic{definition}}

\newenvironment{theo}{\refstepcounter{definition}
\bigskip\noindent{\it Theorem.}\it }{\medskip}

\newcommand{\mar}[1]{}

\begin{document}
\hbox{}

{\parindent=0pt

{\large\bf The quasi-periodic stability condition (the KAM theorem) for
partially-integrable systems}
\bigskip

{\sc G.SARDANASHVILY}
\medskip

\begin{small}

Department of Theoretical Physics, Physics Faculty, Moscow State
University, 117234 Moscow, Russia. 

E-mail: sard@grav.phys.msu.su

URL: http://webcenter.ru/$\sim$sardan/
\bigskip

{\bf Abstract.}
Written with respect to an appropriate Poisson structure, a partially
integrable Hamiltonian system is viewed as a completely integrable system
with parameters. Then, the theorem on quasi-periodic stability in Ref. [1]
(the KAM theorem) can be applied to this system.
\end{small}
}

\bigskip
\bigskip

A Hamiltonian system on a $2n$-dimensional
symplectic manifold $(Z,\Om)$ is said to be 
a partially integrable system (henceforth PIS) 
if it admits $1\leq k\leq n$
integrals of motion $H_i$, including a Hamiltonian,
which are in involution and independent almost
everywhere on $Z$.
Let $M$ be its regular connected
compact invariant manifold. Under the conditions of 
Nekhoroshev's theorem \cite{gaeta,nekh,jpa03},
there exists an
open neighbourhood of $M$ which is a trivial composite bundle
\mar{zz10}\beq
\pi:U=V\times W\times T^k\to V\times W\to V \label{zz10}
\eeq
over domains $W\subset \Bbb R^{2(n-k)}$ and 
$V\subset \Bbb R^k$. It
is provided with the partial
action-angle coordinates
$(I_i,z^A, \f^i)$, $i=1,\ldots,k$, $A=1,\ldots,2(n-k)$, such that
the symplectic form $\Om$ on $U$ reads
\mar{dd26}\beq
\Om= dI_i\w d\f^i +\Om_{AB}(I_j,z^C) dz^A\w dz^B +\Om_A^i(I_j,z^C)
dI_i\w dz^A, 
\label{dd26}
\eeq
and integrals of motion $H_i$ depend only on the action
coordinates $I_j$.

Note that
one can always restrict $U$ to a Darboux coordinate
chart  provided with
coordinates $(I_i,p_s,q^s;\vf^i)$ such that the symplectic form $\Om$ 
(\ref{dd26})
takes the canonical form
\be
\Om= dI_i\w d\vf^i + dp_s\w dq^s
\ee
\cite{fasso,epr2}.
Then, the PIS $\{H_i\}$ on this
chart can be extended to a completely
integrable system, e.g., $\{H_i,p_s\}$, but its invariant manifolds fail
to be compact. Therefore, this is not the case of
the KAM theorem.

Let $\cH(I_j)$ be a Hamiltonian of a PIS in question on $U$ (\ref{zz10}).
Its Hamiltonian vector field
\be
\xi_\cH=\dr^i\cH(I_j)\dr_i 
\ee
with respect to the symplectic form $\Om$ (\ref{dd26})
yields the Hamilton equation
\mar{k65}\beq
\dot I_i=0, \qquad \dot z^A=0, \qquad \dot\f^i=\dr^i\cH(I_j) \label{k65}
\eeq
on $U$. Let us consider a perturbation
\mar{k12'}\beq
\cH'=\cH+\cH_1(I_j,z^A,\f^j) \label{k12'}
\eeq
(see \cite{epr1} for the case of $\cH_1$ independent of $z^A$).
A problem is that the Hamiltonian vector field of the perturbed Hamiltonian
(\ref{k12'}) with respect to the symplectic form $\Om$
(\ref{dd26}) leads to the Hamilton equation $\dot z^A\neq 0$ and,
therefore, no torus (\ref{k65}) persists.

To overcome this difficulty, let us provide the toroidal domain
$U$ (\ref{zz10}) with the Poisson structure of rank $2k$
given by the Poisson bivector field
\mar{k0}\beq
w=\dr^i\w \dr_i. \label{k0}
\eeq
It is readily observed that, relative to $w$,
all integrals of motion of the original PIS $(\Om,\{H_i\})$ remain
in involution and, moreover, they possess the same Hamiltonian vector
fields. In particular, a Hamiltonian $\cH$ with respect to the
Poisson structure (\ref{k0}) also leads to the Hamilton
equation (\ref{k65}). Thus, we can think of the pair $(w,\{H_i\})$
as being a PIS on the Poisson manifold $(U,w)$. The key point is that,
with respect to the Poisson bivector field $w$ (\ref{k0}),
the Hamiltonian vector
field of the perturbed Hamiltonian $\cH'$ (\ref{k12'}) is
\mar{bi51'}\beq
\xi'=\dr^i\cH'\dr_i -\dr_i\cH'\dr^i, \label{bi51'}
\eeq
and the corresponding first order dynamic equation on $U$ reads
\mar{k11}\beq
\dot I_i=-\dr_i\cH'(I_j,z^B,\f^j), \qquad \dot z^A=0,
\qquad \dot \f^i=\dr^i\cH'(I_j,z^B,\f^j). \label{k11}
\eeq
This is a Hamilton equation with respect to the Poisson structure
$w$ (\ref{k0}), but is not so relative to the original symplectic form
$\Om$. Since $\dot z^A=0$ and the toroidal domain $U$ (\ref{zz10})
is a trivial bundle over $W$, one can think of the dynamic equation
(\ref{k11}) as being
a perturbation of the dynamic equation (\ref{k65}) depending on parameters
$z^A$. Furthermore, the Poisson manifold $(U,w)$ is the
product of the symplectic manifold 
\mar{k4}\beq
(V\times T^k,\Om_w), \qquad \Om_w=dI_i\w d\f^i \label{k4}
\eeq
and the Poisson manifold $(W,w=0)$ with the zero Poisson structure.
Therefore, the equation (\ref{k11}) can be seen as a Hamilton equation
on the symplectic manifold $(V\times T^k,\Om_w)$ depending on parameters.
Then, one can apply the conditions of quasi-periodic stability
of symplectic Hamiltonian systems depending on parameters \cite{broer}
to the perturbation (\ref{k11}).

In a more general setting, these conditions can be formulated as follows.
Let $(w,\{H_i\})$, $i=2,\ldots,k$, be a PIS on a
regular Poisson manifold $(Z,w)$ of rank $2k$.
Let $M$ be its regular connected
compact invariant manifold. Under certain conditions \cite{epr2}, there
exists a toroidal
neighbourhood $U$ (\ref{zz10}) of $M$ provided with the partial
action-angle coordinates
$(I_i,z^A, \f^i)$  such that
the Poisson bivector $w$ on $U$ takes the canonical form (\ref{k0}) and
the integrals of motion are independent of the angle coordinates $\f^i$.

We restrict our consideration to the case when all functions are real 
analytic, and they are provided
with the real analytic topology, i.e., the topology of compact convergence
on the complex analytic extensions of these functions. 

We say that a (real analytic) Hamiltonian $\cH(I_j,z^A)$ of a PIS is 
non-degenerate at a point $(\wt I, \wt z)\in V\times W$
if the frequency map
\be
\om:V\times W\ni (I_j,z^A)\mapsto \xi^i(I_j,z^A)\in \Bbb R^k
\ee
is of maximal rank at this point.
Note that $\om(V\times W)\subset \Bbb R^k$ is open
and bounded. Given $\g>0$, let
\be
N_\g=\{ \om\in \Bbb R^k \,:\, |\om^ia_i|\geq \g(\op\sum^k_{j=1}
|a_j|)^{-k-1}, \quad \forall a\in\Bbb Z^k\setminus 0\}
\ee
denote the Cantor set of non-resonant frequencies. The complement
of $N_\g\cap \om(V\times W)$ in $\om(V\times W)$ is dense and
open, but its relative Lebesgue measure tends to zero with $\g\to 0$.

The following result is a reformulation of that in Ref. [1] (Section 5c),
where $P=W$ is a
parameter space and $\si$ is the symplectic form (\ref{k4}) on $V\times T^k$.

\begin{theo} 
Let the Hamiltonian $\cH(I_j,z^A)$ of a PIS in question is non-degenerate 
at a point $(\wt I,\wt z)\in V\times W$, and let 
\mar{k77}\beq
\xi_\cH=\dr^i\cH(I_j,z^A)\dr_i\label{k77}
\eeq
be its Hamiltonian vector field. Then, there exists a 
neighbourhood 
$\cO\subset V\times W$ of $(\wt I,\wt z)$ such that,
for any real analytic locally Hamiltonian vector field
\be
\xi=\xi_i(I_j,z^A,\f^j)\dr^i + \xi^i(I_j,z^A,\f^j)\dr_i
\ee
sufficiently near $\xi_\cH$ (\ref{k77}) in the real 
analytic topology (see [1], pp.59-60), the following holds.
Given the Cantor set $\G_\g=\om^{-1}(N_\g\cap \om(\cO))$, 
there exists a $\xi$-invariant subset of $\cO\times T^k$
which is a 
$C^\infty$-near-identity diffeomorphic (see [1], p.60)
to $\G_\g\times T^k$ and, and on each tori, this diffeomorphism is
an analytic conjugacy from $\xi_\cH$ to $\xi$.
\end{theo}

This is an extension of the KAM theorem to PISs.

\end{document}